# DEPLOYED MODEL OF EXTREMAL SYSTEM OPERATION FOR SOLVING OPTIMAL MANAGEMENT PROBLEMS

I. Lutsenko
PhD, Professor
Department of Electronic Devices
Kremenchuk Mykhailo Ostrohradshyi National University
Pervomaiskaya str., 20, Kremenchuk, Ukraine, 39600
E-mail: delo-do@i.ua

*Розроблено методику побудови моделі розгорнутої системної операції, яка враховує рух кількісних параметрів значущих вхідних і вихідних продуктів операції, знос технологічного механізму, у вигляді одного з вхідних продуктів, і вартісні (експертні оцінки) цих продуктів у часі. Інформація про процес моделі операції доступна в будь-який момент часу*

*Ключові слова: дослідження операцій, модель операції, математичні моделі дослідження операцій, оптимальне управління*

*Разработана методика построения модели развернутой системной операции, которая учитывает движение количественных параметров значимых входных и выходных продуктов операции, износ технологического механизма, в виде одного из входных продуктов, и стоимостные (экспертные оценки) этих продуктов во времени. Информация о процессе модели операции доступна в любой момент времени*

*Ключевые слова: исследование операций, модель операции, математические модели исследования операций, оптимальное управление*

## 1. Introduction

For solving optimal management problems, it is necessary to satisfy a number of conditions in relation to the managed system. The first condition is the possibility of independent management of feeding of raw materials and energy products to the input of the system under study, while ensuring the given quality of a primary useful product at the output of this system. The second condition is the possibility to take into account the wear of the basic technological mechanism in general. The third condition is the possibility of cost or expert estimation of input and output products of the system operation.

This is caused by the need to maximize the desired product of the system in the process of manufacture of the main useful product. The desired product in managed systems is the value added (price), but only extremal systems ensure a purposeful change of the desired product, and optimal systems ensure its maximization.

The above conditions are met in systems with batch feeding of raw materials, at the input of which there are raw material feeding system/systems, and at the output – output technological product buffering system/systems.

## 2. The analysis of literature data and formulation of the problem

Analysis of available publications shows that the continued interest is drawn to the research topic of operations. This indirectly indicates that the central questions remain open.

Conclusions of analytical publication on operations research show that there is a steady trend in demand for works with "rigorous use of empirical data, quantitative analysis and theoretical modeling" [1].

At the same time, it must be recognized that the main efforts are still focused on modeling physical processes of transformation [2–4].

Works, related to inventory management are no exception [5]. Despite the fact that this model is also used for trading systems management, research is focused on the quantitative parameters of products and practically do not affect the cost parameters that directly impact the efficiency of the operation.

Marketers, whose attention is directed at the system process management efficiency [6], which stimulates the development of a generalized model of the operation show great interest in the operations research. However, the lack of a reference point as the efficiency criterion leads to the need to be guided by a set of indicators in researches [7].

Constructing a mathematical model of the process under study is the first and fundamental step towards the operations research. Nevertheless, the analysis of scientific publications indicates that authors, as a rule, without the detailed development of the operation model immediately switch to research methods of private scientific or practical problems that have the ability to involve the developed mathematical apparatus.

For example, in the known work of Hemdi A. Taha [8] statement begins not with construction of a generalized model of the operation, but an introduction to linear programming, a description of the simplex method, etc.

The above analysis can be summarized by the words of Peter Drucker from his legendary work of 1964: "What is a transaction? Above all, how does one decide which of the many transactions within a business is the transaction that is representative of the actual cost structure? There is no set answer" [9].

### 3. The purpose and objectives of the study

For more accurate identification of system operations, it is necessary to develop a mathematical model that links together all the relevant quantitative parameters of technological products and their corresponding expert estimates in time.

This presupposes solving the following tasks:
– the construction of the registration model of the operation;
– the assessment of the model limitation and development of the conceptual model of the extremal managed system with a possibility to determine the wear of technological mechanism and reduce the registration signals of quantitative parameters for technological products to comparable values;
– the construction of the reduced registration model of the operation to comparable values;
– the construction of a generalized, deployed model of the operation;
– the comparative analysis of possibilities of registration and deployed models of operations.

### 4. Development of the deployed model of the system operation

As a result of the interaction of the technological product feeding systems, system under study and output technological product receiver systems, with respect to the system under study, the displacement of technological and information products is formed. Registration of the quantitative parameters of exchange products of these systems allows to create a basic registration system model of the operation.

As an example, let us consider a system operation of batch heating of the liquid using an electric heater. The choice of such a system as a demonstration model is caused by a good knowledge of the model of the heating process, its inertia (there is no need in parameters stabilization nodes) and popularity of wear model of electric heater.

Fig. 1 shows the structure of the managed system, which ensures fulfilling the first condition of optimal management.

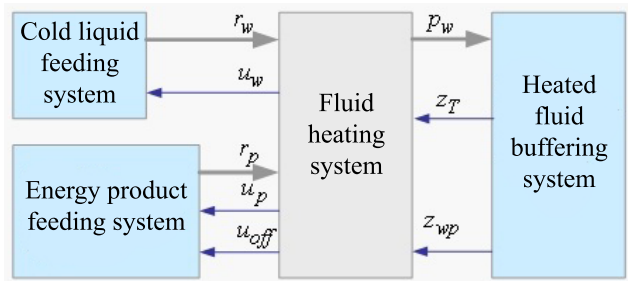

Fig. 1. Structure of the managed fluid heating system with the batch feeding of cold fluid and management of its heating to a predetermined temperature

Following the injection of the setpoint signal $z_{wp}$, associated with the need to replenish the heated fluid buffering system, heating system generates a control signal $u_w$ to feed a certain amount of cold liquid $r_w$. Completion of the process of cold fluid feeding leads to the formation of a control signal $u_p$, which provides the feeding of the energy product $r_p$ with a predetermined intensity. This naturally leads to an increase in temperature of the heated fluid. When the fluid temperature reaches a predetermined value $z_T$, a control signal $u_{off}$ is generated to shut down the feeding of the energy product, and the heated fluid is passed to the buffering system.

A conceptual view of registration signals of quantitative parameters of technological products relative to the heating system is shown in Fig. 2. The model of this type is defined as the registration model of the operation.

Studies show that in energy-consuming technological processes with sufficiently high efficiency, energy consumption decreases with an increase in speed of the technological operation. The operation of fluid heating using the electric heater is not an exception.

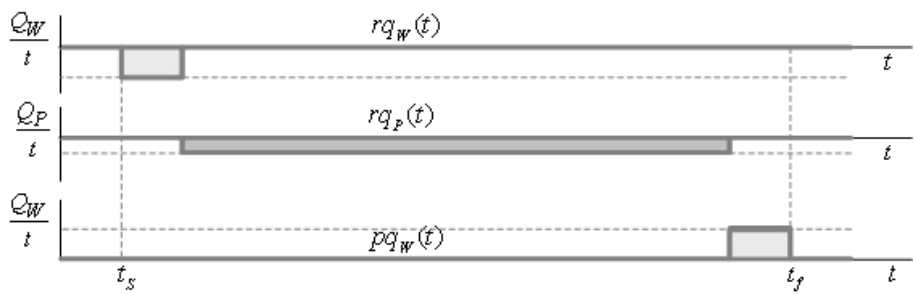

Fig. 2. Timing diagrams of registration signals of technological products at the input and the output of the heating system: $rq_w(t)$ is a registration signal of cold fluid flow rate; $rq_p(t)$ is a registration signal of the energy product consumption; $rq_w(t)$ is a registration signal of the heated liquid flow rate

Let us define the start point of injection of setpoint signal $z_{wp}$ as the start point of the system operation of fluid heating ( $t_s$ ), and the end point of feeding of the heated fluid $p_w$ to the buffering system as a point of physical completion of the system operation of heating ( $t_f$ ). Then the volume of cold fluid $RQ_w$, fed to the input of the heating system is determined by an expression $RQ_w = \int_{t_s}^{t_f} rq_w(t)dt$, where $rq_w(t)$ is a registration signal of movement of cold fluid at the output of the heating system. Similarly, $RQ_p = \int_{t_s}^{t_f} rq_p(t)dt$; $PQ_w = \int_{t_s}^{t_f} pq_w(t)dt$, where $RQ_p$ is the consumption of the energy product during the system operation; $PQ_w$ is the volume of the transferred heated fluid.

Based on the integrated natural indicators of products of the system operation, it is possible to construct the dependences of their change from a given feed rate of the energy product $RQ_w = f(U_P[n])$, $RQ_p = f(U_P[n])$, $RQ_m = f(U_P[n])$ and $PQ_w = f(U_P[n])$.

Studies of the heating process were carried out by its modeling in an EFFLI environment [10]. The results are shown in Fig. 3.

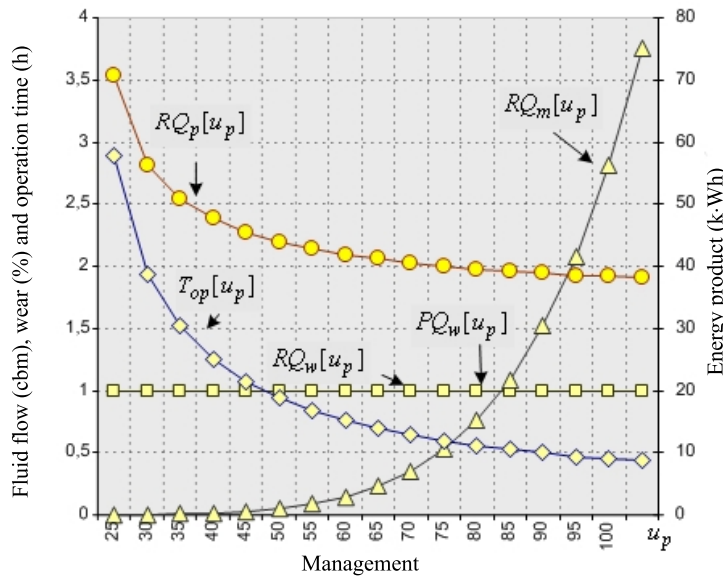

Fig. 3. Integrated flow rates of technological products for system operations of the fluid heating in EFFLI environment, depending on the control $u_p$

Even without a thorough analysis, it is clear that if we do not take into account the wear factor of heating mechanism, feeding of energy product is beneficial to increase to a maximum since the operational energy consumption is reduced along with the time of the system operation.

Studies also show [11, 12] that an increase in the energy product feed rate leads to a disproportional increase in the heating mechanism wear. This means that optimal systems in general should take into account the wear of technological mechanism in constructing the model of the operation of extremal system.

Expression, connecting the energy product consumption and the heating mechanism wear appears as [11] $T = T_н k_u^{-\alpha}$, where $T$ is the life of the heater; $T_н$ is the number of operation hours at rated power; $k_u$ is the ratio of actual output power to its nominal value; $\alpha$ is an indicator, depending on the type of heater. Knowing the power consumption and using this expression, it is possible to calculate the wear of the heating mechanism $RQ_m$ for the particular operation.

Also from the above diagrams it is clear that making a decision on selecting the optimal management is not possible without reducing the quantitative parameters of technological products to comparable values.

Thus, in addition to the introduction of a conditional transmission circuit of more worn equipment, after forming one more system operation, it is also necessary to transmit information regarding the cost estimates of the input and output products to the system input (see Fig. 4).

With this approach, it is possible to reduce the natural quantitative parameters of products of the system operation to comparable values: $re_w(t) = rs_w \times rq_w(t)$; $re_p(t) = rs_p \times rq_p(t)$; $re_m(t) = rs_m \times rq_m(t)$; $pe_w(t) = ps_w \times pq_w(t)$, where $rs_w$ is the cost estimate of a volume unit of cold fluid; $rs_p$ is cost estimate (kWh) of electricity; $rs_m$ is the cost estimate of the wear unit of the heater; $ps_w$ is the cost estimate of the unit of the heated fluid (cbm).

In general, for the i th input product of the operation and the j th output product of the operation we obtain $re_i(t) = rs_i \times rq_i(t)$; $pe_j(t) = ps_j \times pq_j(t)$.

Now we can construct a generalized, reduced, registration model of the operation if we combine streams $re_i(t)$ into a single stream $re(t)$ by the input, and streams $pe_j(t)$ into a single stream $pe(t)$ by the output. Since mathematical notation of summing the functions in time is not available, we introduce the notation $\sum_i f_i(t)dt$ for this operation.

Then $re(t) = \sum_i re_i(t)$; $pe(t) = \sum_j pe_i(t)$.

We define the model of the operation in the form of functions $re(t)$ and $pe(t)$ as the reduced registration model of the operation.

If we replace the function $re_i(t)$ with the impulse function $re_i^*(t)$, where the magnitude of the impulse is determined by the expression $RE_i = rs_i \int_{t_s}^{t_f} rq_i(t)dt$, and the location of the impulse corresponds to the starting point of the operation $pe_j^*(t)$, and if we replace the function $pe_j(t)$ with the impulse function $pe_j^*(t)$, where the magnitude of the impulse is determined by the expression $PE_j = ps_j \int_{t_s}^{t_f} pq_j(t)dt$, we will be able to further simplify the model of the type $re^*(t)$ and $re^*(t)$.

We define the model of the operation in the form of functions $re^*(t)$ and $re^*(t)$ as the registration model of the reduced simplified operation.

Now, the function $re(t)$ can be represented as an impulse signal with amplitude $RE$ and positioned at a point, corresponding to the start point of the system operation $t_S$. The function $pe(t)$ can be represented as an impulse signal with amplitude $PE$ and placed in a point, corresponding to the point of physical completion of the system operation $t_f$ (Fig. 5).

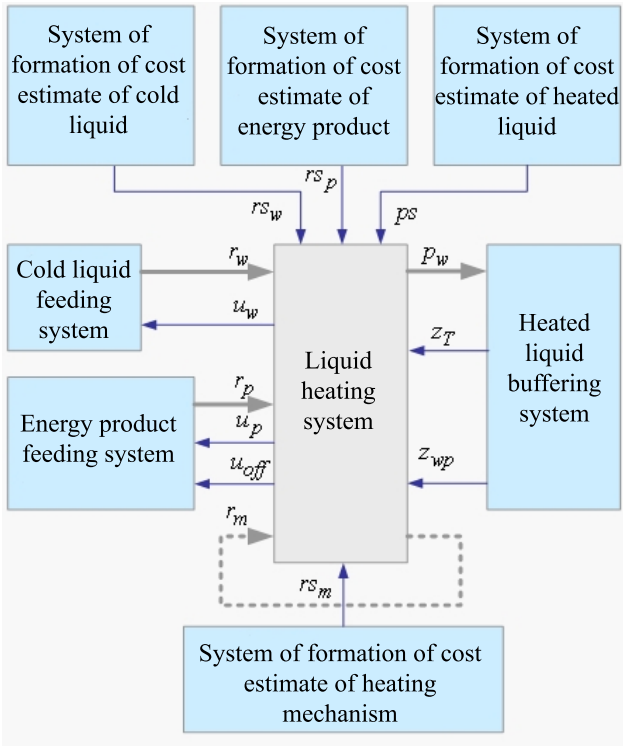

Fig. 4. Block diagram of the extremal heating system with a batch feeding of cold fluid, its heating management to a predetermined temperature, taking into account the wear of the heating mechanism and the cost estimates of technological products

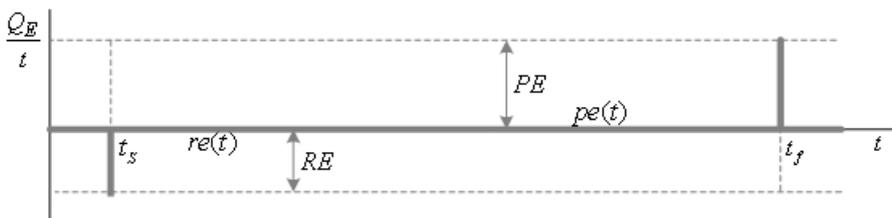

Fig. 5. Registration model of the reduced, simplified operation

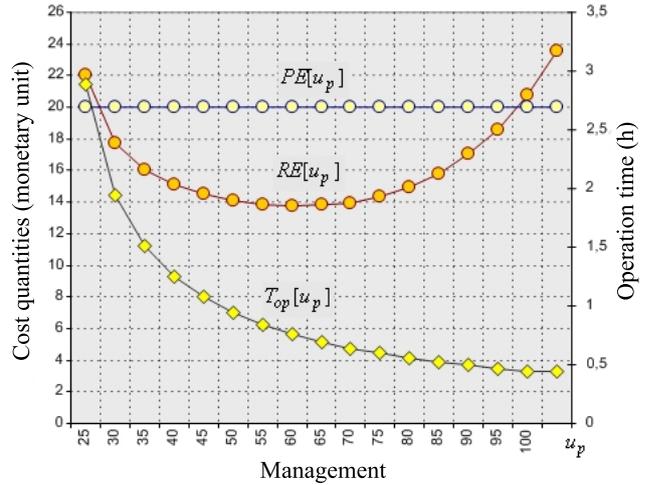

Fig. 6. Comparison of system operations, presented in the form of triple $\theta = (PE, RE, T_{OP})$

Let us turn to our demo system.

In accordance with the foregoing, the volume of cold fluid $RQ_w$, fed to the input of the heating system in comparable values is determined by the expression $RE_w = rs_w \int_{t_s}^{t_f} rq_w(t)dt$, where $RE_w$ is the cost of consumed cold fluid. Similarly $RE_p = rs_p \int_{t_s}^{t_f} rq_p(t)dt$; $RE_m = rs_m \int_{t_s}^{t_f} rq_m(t)dt$; $PE_w = ps_w \int_{t_s}^{t_f} pq_w(t)dt$; where $PE_p$ is the cost of electricity consumed; $RE_m$ is the cost of wear of the heater; $PE_w$ is the cost of the heated fluid.

Creating the reduced, simplified model of the operation allows to define it as a tuple $\theta = (PE, RE, T_{OP})$, and the set of such operations graphically (Fig. 6), where a large number of system operations, each corresponding to a specific feeding of energy product is presented. Here $RE = RE_w + RE_p + RE_m$, $PE = PE_w$.

The study of the data obtained (Fig. 6) shows that the minimum cost management mode ($u_p$=60) is not likely to be the most effective. Increase in the power ($u_p$=65, $u_p$=70) results in a slight increase in costs and a significant reduction in the heating operation time. Therefore, to select the most effective management mode, a more "subtle" tool than triple $\theta = (PE, RE, T_{OP})$ is needed.

Today science knows one recognized relative indicator that identifies the operations based on total cost estimates by the input and output of the system under study – return. If we transfer the ideology of return to our cybernetic model of the operation, where there are no actually operating cash flows, we obtain the following expression for calculating the conditional return (R) $R = (PE - RE) / RE$.

If we try to identify two operations on the Fig.6 with control $u_p$=50 and $u_p$=70, we obtain the same results since their values PE and RE are equal $R = (20 - 14) / 14 = 0,428$. Thus, it is obvious that the operation with control $u_p$=70 is more efficient because the time of this operation is less ceteris paribus.

Therefore, despite the fact that the simplified model of the reduced registration operation contains all the necessary information for the analysis of the studied process, the direct use of the information contained therein does not allow to identify the operation. For example, a performance indicator has not yet been developed on the basis of the registration model of the operation.

Let us consider the reason for this state of affairs.

The fact is that registration signals represent the measured parameter of the studied product with respect to the registrator (sensor, credit/debit document of the bar code scanner sensor, and so on). In this case, registration signals "convey" nothing about the internal state of the process.

For example, we investigate a simplified model of the operation of displacement of a product, linearly distributed in space. We neglect energy costs, associated with the displacement process.

Sensors, registering quantitative parameter of the product, displaced in time (Fig. 7) are installed at the inlet and outlet of technological process of displacement.

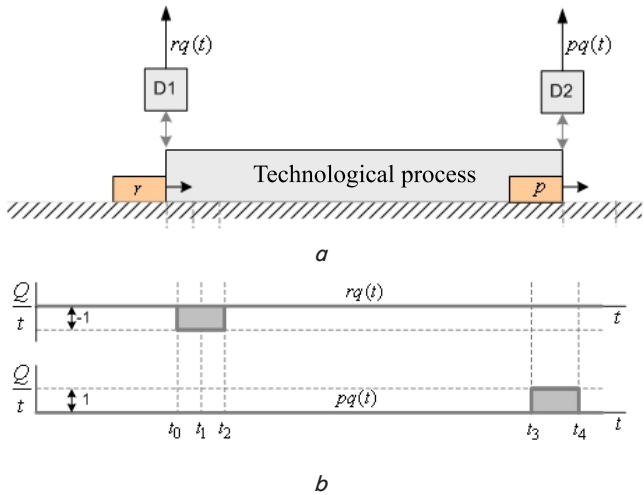

Fig. 7. Displacement process model: *a* — object model; *b* — registration model

As a result of the operation of sensors and display of the results, we obtain registration model of the operation in the form of signals $rq(t)$ and $pq(t)$. This is how the sensors "see" the displacement process. And what does the observer of internal processes see?

In order to answer this question, let us conduct a thought experiment, having "instructed" the observer A to continuously record in time accounting parameter of the visible part of the object under observation.

At time $t_0$ observer does not see the displacement product. After time $t_0$ he sees increasing sizes of object. At moment $t_1$, the observer sees the half of the object, at moment $t_2$ he observes the entire object. Moment $t_3$ is characterized by the fact that it is the last time, when the observer sees the entire object.

The picture that would be displayed by the observer will have a form of the function $icq(t)$ (Fig. 7). Let us define this function as a flow of tight resources. Each point of this flow shows what amount of the displaced product in a given time is bound by the technological process of the operation. The sensor signals do not display such information directly. They are at the input and output of processes of the operation.

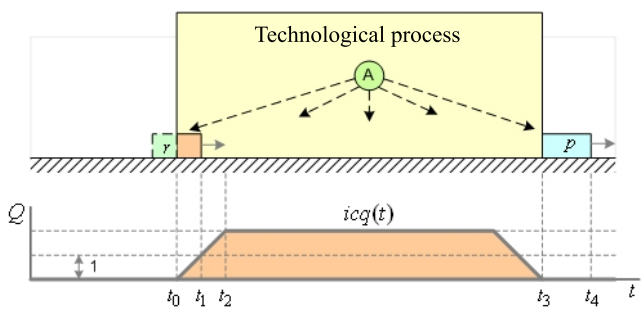

Fig. 8. The principle of constructing the flow of tight resources of the operation from the internal observer's side

In order to get a picture of the internal process, it is necessary to pre-process, integrate registration signals, and then find the difference of the resulting functions.

Let us define the integral of the function $rq(t)$ as $irq(t)$, and the integral of the function $pq(t)$ as $ipq(t)$. Then the sum of the functions $icq(t) = irq(t) + ipq(t)$ will give us a picture that shows the inner processes of the operation (Fig. 9).

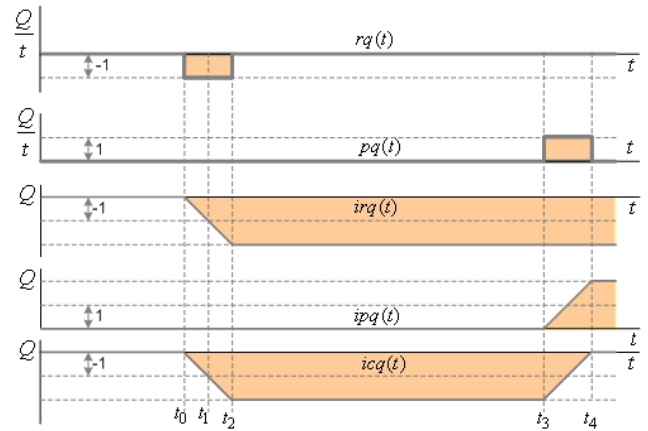

Fig. 9. The deployed quantitative model of the displacement operation of the operation product

These processes are examined from the external observer's side, so the function $icq(t)$ is located in the negative half-plane.

We use deployment technology in relation to the registration model of the simplified, reduced operation (Fig. 10).

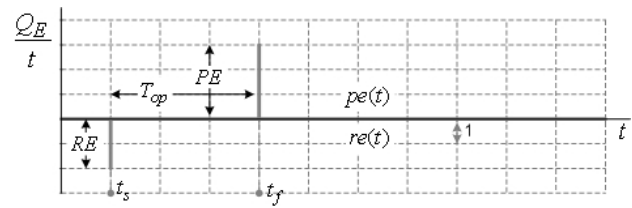

Fig. 10. Registration model of the simplified, reduced operation

To do this, we deploy the model in Fig. 10, having integrated functions $re(t)$ and $pe(t)$: $ire(t) = \int_{t_0}^{t} re(t)dt$; $ipe(t) = \int_{t_0}^{t} pe(t)dt$ (Fig. 11).

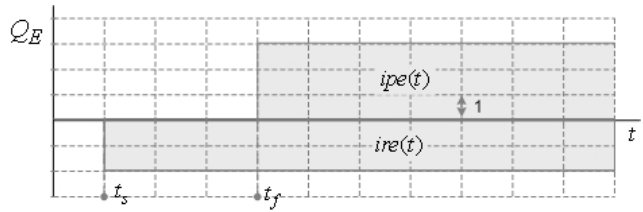

Fig. 11. System model of the operation in the form of reduced combined and deployed registration signals by input and output

Now we can construct the deployed model of the simple reduced operation as a function $ice(t)$ (Fig. 12) using the expression

$$ice(t) = ire(t) + ipe(t) = \int_{t_0}^{t} re(t)dt + \int_{t_0}^{t} pe(t)dt.$$

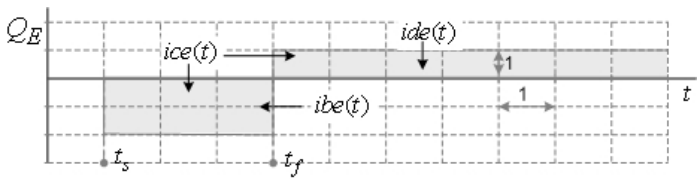

Fig. 12. The deployed model of the reduced, simplified system operation

The function $ice(t)$ can be represented as two functions, the function in the form of a closed flow of tight resources and function in the form of an open flow $ide(t)$. Building the deployed model of the reduced operation based on the original functions with the registration of the time-phased parameters of technological products presents no practical difficulty.

In conclusion, let us consider the possibilities of identifying the operations with equal pairs $PE$ and pairs $RE$ and different time of the operation, but using the new opportunities of deployed model of the operation. In this case, in order to better visualize the process we use the model, shown in Fig. 12.

Let the first operation $PE = 3$, $RE = 2$ and $T_{op} = 3$. Then, for the index, which we denote as ($F$), having integrated the target flow $ide(t)$ on the unit interval and closed flow $ibe(t)$ on the interval $[t_s, t_f]$ we obtain

$$F = \int_{t_f}^{t_f+1} ide(t)dt \Big/ \int_{t_s}^{t_f} ibe(t)dt = 1/6 = 0{,}167.$$

For the operation, in which $T_{op} = 2$, the result is already otherwise $F = 1/4 = 0{,}25$. Since less time of the operation, ceteris paribus indicates a more efficient operation, we see that the index $F$ opens up more possibilities than traditional index $R$.

Thus, the deployed model allows to take into account the "time" parameter when identifying operations.

## 5. Conclusions

On an example of the system operation of fluid heating, process of constructing a basic mathematical model of the operation, based on the registration of the sensor signals, which represent the quantitative parameters of the input and output technological products and taking into account wear of technological mechanism was shown. The concept of the base registration model of the operation was introduced.

Reduced, simplified registration mathematical model of the system operation, obtained by scaling the registration signals of the basic model and determining the sum of the reduced, integrated flow rates by the input and output of the system under study was constructed. The concept of the reduced registration model of the operation was introduced.

The main limitation of the registration model of the operation, associated with the features of primary data collection was illustrated, and a way out of this situation by deploying a registration model to get access to data of operations at each point of the technological process of the operation was shown. The concept of the deployed model of the operation was introduced.

The developed approach to constructing the mathematical model of the reduced deployed operation allows more accurate identification of operations in comparison with known methods. Capabilities of the deployed model of the operation allow to use it to obtain new scientific results in the field of operations research.